\documentclass[12pt]{amsart}
\usepackage{amsmath}
\usepackage{amsthm}
\usepackage{amsfonts}
\usepackage{amscd}
\usepackage{amssymb}
\usepackage{mathrsfs}
\usepackage{graphicx}
\usepackage{psfrag}

\def\ca{{\mathcal A}}

\def\cj{{\mathcal J}}

\def\cam{{\mathcal M}}
\def\cn{{\mathcal N}}

\def\cs{{\mathcal S}}
\def\ct{{\mathcal T}}

\def\bb{{\mathbb B}}
\def\bc{{\mathbb C}}

\def\bn{{\mathbb N}}
\def\br{{\mathbb R}}
\def\bz{{\mathbb Z}}

\def\d{\delta}        
\def\eps{\varepsilon}

\newtheorem{Thm}{Theorem}[section]

\newtheorem{Dfn}[Thm]{Definition}

\newtheorem{exmp}[Thm]{Example}

\theoremstyle{remark}
\newtheorem{rem}[Thm]{Remark} 
\newtheorem{ack}{Acknowledgement} 

\begin{document}

 \title{ Sums of two dimensional spectral triples }

 \author{ Erik Christensen and Cristina Ivan}

 \address{  
  Department of Mathematics, University of Copenhagen,   Denmark \\ 
  Department of Mathematics, University of Hannover,  Germany }

 \email{  echris@math.ku.dk, antonescu@math.uni-hannover.de}
 \date{\today}


 \begin{abstract}
   We study countable sums of two dimensional modules for  the
   continuous complex functions on a compact
   metric space and show that it is possible to
  construct a spectral triple 
   which gives the original metric back. This spectral triple will be
   finitely summable for any positive parameter. We also construct a
   sum of two dimensional modules which reflects some  aspects of the
   topological dimensions of the compact metric space, but this will
   only give the metric back approximately. At the end we make an
   explicit computation of the last module for the unit interval in
   $\br.$  The metric is recovered exactly, the Dixmier trace 
   induces a multiple of the Lebesgue integral  and 
   the growth of the number of eigenvalues $N(\Lambda)$  bounded by $\Lambda$
   behaves, such that $N(\Lambda)/\Lambda$ is bounded, but without
   limit for $\Lambda \to \infty .$
\end{abstract}

\maketitle

\section{Introduction}
\noindent
Gelfand's fundamental theorem on Abelian C*-algebras shows that the
study of unital commutative C*-algebras is the same as the study of
compact, Hausdorff, topological spaces. It has been clear for many
years 
that the non commutative C*-algebras have many properties in common
with their Abelian relatives, for instance the theory of measure and
integration was already well developed in the non commutative setting
by Murray and von Neumann. The basic difference is of course that for
a unital Abelian C*-algebras there is a space, the spectrum,  which can be
investigated via lots of different mathematical theories. Quite a few
theories can on the other hand be expressed in terms of  certain
subalgebras of the algebra of continuous functions on the spectrum. As
 examples, one may  think of differentiability and
  Lipschitz continuity. Also other
structures like K-theory are expressible in terms certain subalgebras of
the continuous functions on the spectrum and matrix algebras over
such algebras. Especially Alain Connes has tried to express
geometrical structures this way in order to be able to extend the
classical differential geometry to the non commutative world. In the
classical setting  a differential geometric
structure can be viewed as something -  extra - 
added to the compact topological space
in question. In the  non commutative case one should expect that 
there is a C*-algebra in the background and that the geometrical
features such as smooth functions is something extra.
 In Connes' work \cite{Co1} this has
been synthesized in the concept called a spectral triple, which is
defined as

\begin{Dfn} \label{spectrip}
Let $\ca$ be a unital C*-algebra, $H$ a Hilbert space which carries a
unital representation $\pi$ of $\ca$ and $D$ an unbounded self-adjoint
operator on $H$. The set $(\ca, H, D) $ is called a spectral triple if
\begin{itemize}
\item[(i)] the set $\{a \in \ca\,|\, \|\,[D, \pi(a)]\, \| < \infty
  \,\}$ is a dense subset of $\ca$,
\item[(ii)] the operator $(I+D^2)^{-1}$ is compact.
\end{itemize} 
If {\rm tr}$((I+D^2)^{-p/2}) < \infty $ for some positive $p$ then the
spectral triple is said to be $p$-summable, or just finitely
summable. 
\end{Dfn}
\noindent
The set  $\{a \in \ca\,|\, \|\,[D, \pi(a)]\, \| < \infty
  \,\}$ may be thought of as the set of continuous functions which
  have essentially bounded derivatives, but the concept of a spectral
  triple is much more than a way to express an analogy. One of the
  striking results by Connes is the theorem which
  shows that the geodesic distance on a compact spin Riemannian
  manifold can be computed via a spectral triple. This setup is
  described in \cite{Co2} Chapter VI. The expression which
  gives the geodesic distance back on a manifold can easily be
  extended to define a metric on the set of regular Borel probability 
measures on the manifold. On the other hand such probability measures
are exactly what is called states on a C*-algebra so we get in this
way a metric on the state space of a C*-algebra from a spectral
triple, \cite{Co2} VI.1, p. 544:

\begin{Dfn} \label{disttrip}
Let  $ST=(\ca, H, D) $ be a spectral triple then the induced metric
$d_{ST}$ on the state space $\cs(\ca)$ is 
defined by:
\begin{displaymath}
\forall \phi, \psi \in \cs(\ca): \quad d_{ST}(\phi, \psi)\, := \,
\sup\,\{\, | \, \phi(a) - \psi(a)\,|\,|\, \|\, [D, \pi(a)]\, \|\, \leq
\, 1\, \}.
\end{displaymath}
\end{Dfn}
\noindent
It should be remarked that such a  metric is non-standard in the way
that it may take the value $\infty$ on certain pairs of states, but it
is on the other hand fascinating, that the non commutative world keeps
track of the metric on a space, which no longer exists in it's
classical form.

\medskip
\noindent
In the present paper we study especially  {\em the metric of
  compact metric  spaces from a non
commutative point of view.  } This kind of investigation has already
been performed by Marc 
Rieffel. In several papers, among which we only cite a few
\cite{Ri1, Ri3}, Rieffel studies 
 a general compact metric  space by considering a subalgebra of the
algebra of continuous functions on the space and a seminorm on this
algebra, rather than studying the space and a metric on the space. 
He also  showed how the metric of any compact metric space 
can be recovered from a {\em  Dirac} operator which is a sort of a
differential operator. This construction does not give a spectral
triple because the Dirac operator will not have compact resolvent.
 
\medskip
\noindent
Connes has associated the concept called a  Dixmier trace to a
spectral triple, and shown that in the case of a compact spin
Riemannian manifold this is a multiple of the volume form \cite{Co2},
Formula 2, p. 545. More than
10 years ago Michel L. Lapidus realized that Connes' spectral triples
and associated Dixmier traces are important tools in his ongoing 
search for possible extensions of classical results by Weyl on the
asymptotic growth of the eigenvalues of the Laplacian operator on a
bounded open subset of $\br^d$ with a {\em nice }  boundary. Lapidus
wanted to extend these results to bounded open subsets of $\br^d$ with
fractal boundaries \cite{La1} and even further to fractal subsets of
$\br^d$ \cite{KL, La2, LM, LP}. In these papers the study of the
geometry is
linked to  an investigation of the Laplacian and its spectral
distribution. This scope is a little bit different from ours since we are
looking at possible Dirac operators, whose squares are supposed to be 
Laplacians. On the other hand a spectral triple,
as constructed in this article, is of a
discrete nature, since the associated 
representation of the algebra of continuous
functions takes place inside an algebra of type $\ell^{\infty}.$ 
For a classical differential operator the associated
representation of the function algebra is usually as an ultraweakly
dense subalgebra of an $L^{\infty}-$algebra. 
Anyway the results of Weyl can  make sense in both settings, since
both settings come with Dirac operators with  discrete spectra. 

\smallskip
\noindent
In the papers \cite{La2, LM, LP} Lapidus, Maier and
Pomerance   establish a connection
between the Riemann Hypothesis and the Weyl-Berry Conjecture for a
fractal string. The relations between Riemann's zeta function and
fractal geometry  is further developed in the book  \cite{LF}
{\em Fractal geometry and number theory, } by Lapidus and
Frankenhuysen. The possibility to continue this investigation based on
Connes spectral triples and the Dixmier trace lies in the fact that
besides the connection between the Dirac operator and the Laplacian in the
classical case, the zeta  function $\zeta(s)$ is a multiple of 
the trace of the operator $|D|^{-s}$
and the Dixmier trace gives the volume form. The papers \cite{La3,
  La4} contain detailed studies of the ways non commutative geometry
can be  used to investigate important problems  in fractal geometry.    

\smallskip
\noindent
In this paper we are studying abstract compact metric spaces and
included here are compact fractal subsets of $\br^d$, and  we hope that
 the spectral triple we are proposing in
Theorem 3.2 may turn out to be usable in Lapidus' project of studying
fractal sets via concepts from non commutative geometry. On the other
hand the spectral triples of this articles Theorem 3.1 and Theorem 3.2
 will give Laplacians of a quite different
nature than the ones studied by Lapidus and his coauthors. Still we
think that there may be some hope that our modules may be usable  in
Lapidus' project. 
In order to demonstrate this point of view we have made an
explicit construction of our module based on Theorem 3.2 
 for the unit interval and shown
that the spectral triple reproduces the usual metric, the Dixmier
trace becomes a multiple of the usual Lebesgue integral and the growth
of the number of eigenvalues  $N(\Lambda)$, numerically bounded by
$\Lambda$,  is such that $N(\Lambda)/\Lambda$ is bounded´, but without
limit for $\Lambda \to \infty$. This very last result shows that, in
general, 
spectral triples are not uniquely determined, even though they may be
constructed by  fairly reasonable algorithms. The Dirac operator of
this example will not satisfy  Weyl's  asymptotic formula 
  in the way that there
exists no constant $c$ such that $N(\Lambda) - c\Lambda $ is of lower
order than $\Lambda.$ The reason for the failure has to be
  sought in the geometric growth of the numerical value of the
eigenvalues of the constructed Dirac operator. The $n'$th eigenvalue,
according to numerical 
size is roughly $2^n$ and its multiplicity is also $ 2^n.$ Hence the
gaps in the spectrum grow exponentially, but so do the
multiplicities too, and that makes $N(\Lambda)/\Lambda $ a bounded
function without a limit for $\Lambda \to \infty.$

\medskip
\noindent
In the book \cite{Co2} Example 2 a), p. 563, Connes introduces
a spectral triple for a set consisting of just two points, and in
\cite{Co3} this construction is applied to a countable family of pairs
of points from the classical {\em middle third } Cantor
set in the unit interval. The direct sum of these two-dimensional
modules and their Dirac operators gives a spectral triple for the
algebra of 
continuous functions on the 
Cantor set, and it is possible to obtain a lot of
exact geometric data
for the Cantor set from this spectral triple. 
This construction has been studied in details and extended to certain 
fractal subsets of $\br^n$ by D. Guido and T. Isola  in  
\cite{GI1, GI2, GI3}. Their constructions are based on the
fact that many fractals in $\br^n$ can be obtained as limits for
iterated function systems of contractions. As in Connes' construction
they can then obtain a good description of the fractal in question by
looking at  a pair of points and the sequence of images of this pair,
generated  by the iterations of the function system. They also have
another construction based on a single point. The pairs  will then be
of the form {\em (child, parent)} 
among the iterated images of the starting point.

\smallskip
\noindent
In this paper we are studying 
an abstract compact metric space and we want to see to which extend it can be
described using a spectral triple which is a countable direct sum of
two-dimensional modules. 
 We show that the metric
can be recovered exactly from such a  spectral triple.
 This construction shows that a metric is relatively easy
to describe via a spectral triple, and that it is difficult to put a
ranking to spectral triples for general compact metric spaces. 
The spectral triple which gives the metric back exactly
is {\em in principle } optimal, because no information regarding the
metric space is lost when we compare {\em the algebra of continuous
functions plus the spectral triple} with  {\em the compact space plus
the metric.  } On the other hand the spectral triple under
consideration is not the {\em right } one because it is summable for
any positive value. In general the spectral triple should somehow
reflect the dimension of the space and the {\em local density } of
the space, as it has been nicely pointed out 
in \cite{GI2}. For this reason we give another construction of a spectral
triple, which is much more complicated and does not give the metric
back exactly, but a metric which is only  Lipschitz equivalent, to the
original one. This spectral triple reflects the local structures of
the space and gives some computable estimates of the upper Minkowski
dimension of the space. In concrete cases it is possible to do much
better as shown by Guido and Isola \cite{GI1, GI2, GI3}. We have tried
our spectral triple on the unit interval and we get the expected
results {\em except  for the growth of the eigenvalues } of the Dirac
operator. The number $N(\Lambda)$ of eigenvalues numerically dominated
by $\Lambda$ is of the order $O(\Lambda)$, but $N(\Lambda)/\Lambda$
has no limit for $\Lambda \to \infty$. It is possible to obtain
similar results for the unit cube in $\br^d$, but the details for the
unit interval are already many.

\begin{ack}
The authors want to thank professor Michel L.  Lapidus from University of
California at Riverside for much good advice and many interesting
discussions on possible applications of the methods from non
commutative geometry to the study of the geometry of fractals.
\end{ack}

\section{Notation and definitions}
\noindent
The paper deals with elementary aspects of metric spaces and uses some
of the language of Connes theory of non commutative geometry 
 to do so. The standard reference to Connes work is his book
 \cite{Co2}. We will also use some basic results on operator algebras,
 and all what we use can be found in the books by Kadison and Ringrose
 \cite{KR},  but many other text
books on operator algebra will also describe the concepts we are
using. 

\medskip
\noindent
In the introduction we referred to Connes' construction of a spectral
triple based on two points \cite{Co2} p. 563, 2. Example a). 
We will use this construction quite a lot
and call such a spectral triple a two-point spectral triple. 
The basic idea is that for a subset
consisting of two different points $\{x,y\}$ of a compact metric space
$(\ct, d)$ 
it is possible to construct a spectral 
triple which can express the distance between these points. It is done
via the following definition.

\begin{Dfn} \label{twopttrip}
Let $(\ct,d)$ denote a compact metric space and  $x, y$ be points in
$\ct$. The two-point spectral triple  $ST_{x,y}\, := \, ($C$(\ct),
H_{x,y}, D_{x,y})$ is defined by
 
\begin{itemize}
\item[(a)] If $ x = y$ then $ H_{x,y} = \,\{0\}.$ 
\item[(b)]  If $ x \neq y \text{ then }    H_{x,y} = \,\bc^2
  \newline \text{ and } \pi_{x,y}: \text{C}(\ct) \to B(H_{x,y}) 
\text{ a representation given by.} $
\begin{itemize}
\item[(i)]
$\text{  for } f \in \text{C}(\ct) \text{ and }
  (\alpha , \beta) \in H_{x,y} \newline \quad \pi_{x,y}(f)(\alpha, \beta) :=\,
  (f(x)\alpha,  f(y)\beta ) $
\item[(ii)] 
$D_{x,y} : H_{x,y} \to \, H_{x,y} \newline
\text{ and } 
D_{x,y}(\alpha, \beta)=
\left ( \frac{\beta}{d(x,y)} ,\frac{\alpha}{d(x,y)} \right ).$
\end{itemize}
\end{itemize}
\end{Dfn} 
\noindent
We remark that  for the metric space
$(\ct, d)$ it follows by an elementary calculation that
the metric $d_{ST_{x,y}}$ which by 
 Definition \ref{disttrip} is induced on $\ct$  will be
given by the following table.
\begin{displaymath}
d_{ST_{x,y}}(u,v) = \begin{cases} 0 \quad \text{ if } \, u = v \\
d(x,y) \quad \text{ if } \{u,v\} = \{x,y\} \\
\infty \quad \text{ if }  u \neq v \text{ and } \{u,v\} \neq  \{x,y\}
\end{cases}
\end{displaymath}

\medskip
\noindent
We will use some ingredients of the theory of dimension of metric
spaces. We will not repeat much of it,  but refer the reader to the book
by Falconer \cite{Fa}.  In order to use this language we will use the
following notation.

\begin{Dfn} \label{ball}
Let $(\ct, d)$ denote  a metric space, then for $r \geq 0 $ and $ t
\in  \ct $ we will
let $\bb(t,r) $ denote the closed ball of radius $r$ centered at $t$.
\end{Dfn}
\noindent
For a compact metric space $(\ct,d)$ and a positive real $r$ the
compactness implies that $\ct$ can be covered by a finite number of
closed balls of radius $r$. Such a covering will be called minimal if
the number of balls is minimal among all finite coverings by closed balls of
radius $r$. The minimal number in such a covering 
is denoted $N_r$ and we remind the reader of the following definition

\begin{Dfn} \label{upmindim}
The upper Minkowski dimension - or upper box counting dimension - is
defined as 
\begin{displaymath}
\overline{\text{dim}_{\text{M}}}(\ct)\, := \, \underset{r\to
  0}{\limsup}\frac{\log(N_r)}{-\log(r)}.
\end{displaymath} 
\end{Dfn}

\noindent
This supremum may not be finite, but if the set $\ct$ is a subset of $\br^d$
and the metric on $\ct$  is the one inherited from $\br^d$ then the
upper Minkowski  dimension is at most $d$.
 
\section{Sums of two-dimensional triples.}  
\noindent
This section contains two examples of spectral triples associated to a
compact metric space. Both triples are countable sums of
two-dimensional modules and both induce metrics for the given
compact topology. The first example gives the metric back exactly and
the second triple depends on  a parameter $\delta > 0$, and it can 
 be constructed such that the induced metric is
within   a $\delta-$distance of the original one. Then why bother with
the second example. The reason is that the first triple is finitely
summable for any positive real number $s$. Hence this spectral triple
contains no information on the dimension of the space. The other
spectral triple reflects dimension properties such as the upper
Minkowski dimension. It is much harder to construct and not so precise
with respect to the metric, but it is probably a more relevant module
for topological investigations. One can then ask what is the
difference. For us it seems that the explanation for this phenomenon
has to be found in the fact that the module $H$, we are using, is a
countable sum of two-dimensional Hilbert spaces. This means that the
operators  in the representation of the continuous functions on the
compact metric space are all in a discrete maximal Abelian von Neumann
subalgebra $\ca$ of $B(H)$. If $D$ is a Dirac operator for the
module, and $N$ is a self-adjoint operator affiliated with $\ca$,
such that $N+D$ has a self-adjoint closure, then for any operator $a
\in \ca$ we have $[D,a] \,= \, [D+N,a]$. Now such an $N$ can be chosen
relatively freely and such that $\left ( I+(D+N)^2 \right ) ^{-s/2}$ is
of trace class for any positive $s$.

\begin{Thm} \label{extripd}
Let $(\ct,d)$ denote a compact metric space, then
 there  exists a sequence   of pairs of
 non-equal points $\{x_n,y_n\} $ from $\ct$ such that the sum of all
 the two-point spectral triples $ST_{x_n,y_n} $ becomes a spectral triple
 for C$(\ct)$. This triple is denoted $ST(d)$, it induces the given
 metric d and it is summable for any positive real $s$.
\end{Thm} 

\begin{proof}
Since $\ct$ is compact and metric it contains a dense sequence
$(t_i)$. We can then let the sequence $\{x_n,y_n\}$ be a numbering of
the set of {\em unordered pairs } 
$\{\{t_i,t_j\} \, | \, i,j \in \bn \text{ and }t_i \neq
t_j\,\}$. For each pair $\{x_n,y_n\}$ we consider the two-point module
$ST_{x_n,y_n}$, but we will modify $D_{x_n,y_n}$ such that new Dirac
  operator, which we call $D_n$ is given by the following 2
  by 2 matrix
\begin{displaymath}
D_n:= \begin{pmatrix} 2^n  &  \frac{1}{d(x_n,y_n)}\\
                                \frac{1}{d(x_n,y_n)}  & -2^n
\end{pmatrix}
\end{displaymath}
The other items are fixed, but in order to ease the notation,
 we will let $H_n$ denote the two-dimensional Hilbert
space $H_{x_n,y_n}$ and $\pi_n$ is the representation $\pi_{x_n,y_n} $
  of C$(\ct) $ on $H_n $. It is easy to see that 
\begin{displaymath}
\forall f \in \text{C}(\ct); \quad [D_n, \pi_n(f)]\, = 
 \begin{pmatrix} 0  &  \frac{f(y_n) - f(x_n)}{d(x_n,y_n)}\\
                                \frac{f(x_n)-f(y_n)}{d(x_n,y_n)}  & 0 
\end{pmatrix}.
\end{displaymath}

\smallskip \noindent
Let us start by showing that there exists a dense set of functions $f$
in C$(\ct)$ such that the commutators $[D,\pi(f)]$ are all bounded and
densely defined. To this end we define 
\begin{displaymath}
\cn \,: =\, \{f \in
\text{C}(\ct) \, | \, \exists y \in \ct \, \forall x \in \ct\,\,  f(x):=
d(x,y) \, \}.
\end{displaymath}
It is a simple consequence of the triangle inequality 
that for any  of the pairs 
$\{x_n,y_n\}$  and for any $f(u)\,:=\, d(u,v)$  in $ \cn$ we have
\begin{displaymath} 
\|\,[D_n,\pi_n(f)]\,\| \, = \,\frac{ |\,d(x_n,v)
  -d(y_n,v)\,|}{d(x_n,y_n)}\, \leq \, 1 
\end{displaymath}
Hence for each $f$ in  $\cn$ we have $\|\, [D,\pi(f)]\,\| \, \leq \, 1$ and
since this operator clearly is defined on the linear span of the
$H_n $ spaces, it is densely defined too.  Since for any two points
$u \neq  v $ we have $d(u,v) > 0 = d(v,v)$, we see that the algebra    
generated by the functions in $\cn$ and the constant function $I$
 separates the points in $\ct$, and
hence by Stone-Weierstrass' theorem it is uniformly dense. 
On the other hand it follows from the
derivation property of a commutator that all elements in this algebra have
bounded commutators with $D$, so the first condition for $ST(d) $
to be a spectral triple is fulfilled. Before we start to prove that
the second condition holds, we will collect some  of the results above for
later use in  the following statement.

\begin{equation} \label{msepst}
\text{Let } \cam \,:= \{g \in \text{C}(\ct)\, |\, \|\, [D, \pi(g)]\,\|
\leq 1 \, \} \text{ then } \cam \text{ separates the states}.
\end{equation}

\noindent
Let us then return to the definitions of the $D_n$ operators and of
$D$. First we compute the eigenvalues of each $D_n$, and we get the
set below.

 \begin{displaymath}
\sigma(D_n)\, = \, \left \{ \, -\sqrt{2^{2n}+
  d(x_n,y_n)^{-2}}, \quad \sqrt{2^{2n}+ d(x_n,y_n)^{-2}}\, \right \} 
\end{displaymath}
\noindent
So the two eigenvalues  are both
numerically bigger than $2^n$ and we can conclude that for the
operator $(I + D^2)^{-1}$ there are at most $2n$ eigenvalues of
absolute value bigger than $2^{-n}$. This means that $(I+D^2)^{-1} $
is compact and that for any positive real number $s$ we have 

\begin{displaymath}
\text{tr}\left ( \left ( I+D^2 \right )^{-s/2} \right ) \, < \, 2\sum_{n=1}^{\infty}2^{-ns} \, < \,
\infty.
\end{displaymath}
It then follows that $(\text{C}(\ct), H, D) $ is a spectral triple
which is summable for any positive $s$.

\medskip
\noindent
We will now show that the metric, say $d_{ST(d)}$ induced by $ST(d) $ on
the state space has the property that it agrees with the original
metric on the pure states, i. e. the point measures. We will first
show that for any pair of points $\{x,y\}$ from $\ct$ we have 
$d(x,y) \leq d_{ST(d)}(x,y) $. This type of inequality  is quite
general
 and follows from the
fact that the set of functions $\cn$
is a subset of the set $\cam$. To be more specific, 
let $s, t $ be given points in $\ct$ and
define the function $g(u) \,:=\, d(u,t)$ in $\cn$, then 
\begin{align*}
d_{ST(d)}(s,t) \, =& \, \underset{f \in \cam}{\text{sup}}|f(s)-f(t)|\, \geq
\,\underset{f \in \cn}{\sup}|f(s)-f(t)| \\ \geq&
\,|g(s)-g(t)|\, = \, |d(s,t) - d(t,t)| \, = \, d(s,t).
\end{align*}  

\medskip
\noindent
Let us then show the inequality $d_{ST(d)}(s,t) \leq d(s,t).$
Let again $s,t$ be given from $\ct,$ and let us suppose first that
\begin{displaymath}
d_{ST(d)}(s,t) > \text{ diam}(\ct) \, := \sup\{d(u,v)\, |\, u, v \in \ct\, \}.
 \end{displaymath}
Then there must exist a function $ f \in \cam$ such that $|f(s) -
f(t)| \, > \text{ diam}(\ct)$. Let then $\eps > 0 $ be chosen such
that $|f(s) -
f(t)| - \text{ diam}(\ct) > 3 \eps$, and find a natural number $n$
such that the  pair $(x_n, y_n) $ is so close to the pair $(s,t)$
that $|f(s)-f(x_n)| < \eps$ and $|f(t)-f(y_n)| < \eps$. Since $f \in
\cam$ we have $|f(x_n) - f(y_n)| \leq d(x_n,y_n) \leq \text{
  diam}(\ct)$ and we get the contradiction
\begin{displaymath}
\text{diam}(\ct) + 3 \eps < |f(s) - f(t)| \leq |f(x_n) - f(y_n)| +
2\eps \leq \text{diam}(\ct) + 2 \eps. 
\end{displaymath}
We now know that $d_{ST(d)}(s,t) \leq \text{diam}(\ct) < \infty$, so
for any positive $\eps$ 
there exists a function $f \in \cam$ such that $|f(s) - f(t)| >
d_{ST(d)}(s,t) - \eps/5$. Again, using the continuity of $f$ and the
original metric $d$,
  we find a pair $(x_n, y_n) $ from the
sequence upon which the spectral triple is build, such that $|f(s) -
f(x_n)| < \eps/5$, $|f(t) - f(y_n)| < \eps/5$, $d(s,x_n) < \eps/5$ and
$d(t,y_n) < \eps/5$. Having this we
finally may conclude as follows  
\begin{align*}
d_{ST(d)}(s,t) <& |f(s) - f(t)| + \eps/5 \\  <& |f(x_n) - f(y_n)| +
3\eps/5 \\ \leq &
d(x_n,y_n)+ 3\eps/5\\ <& d(s,t) + \eps.
\end{align*}
So for any pair $(s,t)$ of points from $\ct$ we have $d(s,t) =
d_{ST(d)}(s,t) $. We are nearly done, but we think it is appropriate
to note that the metric $d_{ST(d)} $ really is a metric on the
state space and that it generates the w*-topology on that space. This
can be seen in many ways. We will use the method which Rieffel has
described  in \cite{Ri2} and we have used in \cite{AC1}, Proposition
3.2.
 Consider again the set $\cam $, fix an
element $v \in \ct$ and define a subset $\cam_v$ of $\cam$ by
\begin{equation} \label{camv} 
\cam_v \, := \, \{f \in \cam \, |\, f(v) = 0\,\}
\end{equation}
We know already by (\ref{msepst}) that $\cam$ separates the states of
C$(\ct), $  so according to Rieffel we then just have to show that  $\cam_v$ is
relatively compact in C$(\ct).$ 
This follows from the construction of $\cam_v$,
since by the definition we have  for
any $ f \in \cam_v$ and any pair of
points $\{s,t\} $ from $\ct$ that  $|f(s) - f(t) | \leq
d_{ST(d)}(s,t) = d(s,t)$. This means that the functions in $\cam$ are
equicontinuous. Further for any $t \in \ct$ and any $ f \in \cam_v$ we
have 
$|f(t)| = |f(t) - f(v)| \leq d(t,v) \leq \text{ diam}(\ct)$, so the
set $\cam_v$ is also bounded. By Arzel\`{a}-Ascoli's Theorem we get
that $\cam_v$ is a relatively compact subset of C$(\ct)$ and the
theorem follows.    
\end{proof}

\begin{Thm} \label{apptripd}
Let $(\ct,d)$ denote a compact metric space and let $0 < \delta $ be a
real number. Then there  exists a countable set $\cj$ of pairs of
 non-equal points $\{x,y\} $ from $\ct$ such that the sum of all
 the two-point spectral triples $ST_{x,y} $ over $\{x,y\} \in \cj$
becomes a spectral triple
 for C$(\ct)$. This triple is denoted $ST(\delta)$ and it induces a
 metric  $d_{\delta} $ on $\ct $ such that 
\begin{displaymath}
\forall s, t \in \ct:  \quad d(s,t) \, \leq \, d_{\delta}(s,t) \,
\leq \, (1+\delta)d(s,t).
\end{displaymath}   

\medskip
\noindent
If the upper Minkowski dimension $\overline{\text{dim}_{\text{M}}}(\ct)$
is finite then  the module is finitely summable for any real 
 $s$ such that $s \, >  \, 2\overline{\text{dim}_{\text{M}}}(\ct)$. 
 
\noindent
If the module is summable for some $s > 0$
 and the topological  space $\ct$ is connected, 
then $\overline{\text{dim}_{\text{M}}}(\ct)$ is at most $s$.
\end{Thm} 

 \begin{proof}
The construction can be done in many ways and it may be reasonable to
give a proof which reflects this. In this way the proof becomes
a little less transparent, but more easily applicable to different
metric sets of fractal type. The proof is based on a sequence
of positive numbers $(r_n),  r_n = \theta \rho^{n-1} $ for
some strictly positive real $\theta \text{ and a real } \rho$ such that 
  $0 < \rho < 1$. 
The simplest
argument for the general case 
can probably be obtained when $\theta = \text{diam}(\ct)$, the
diameter of $\ct$,  
and $\rho = 1/2$, but if one
looks at the usual middle-third Cantor set it turns out that $\theta
= 1/2$ and $\rho = 1/3$ will be the most natural choice. On the other
hand it is quite clear that other sequences $(r_n)$ of positive reals
which are finitely summable may be used in certain cases. We have not
been able to find a proof which works in this generality, so the proof
here will be based on sequences of geometric descent.

\smallskip \noindent
Let us then suppose that  some positive reals $\delta, \theta  \text{
  and }  0 < \rho < 1 $  are given.
 We are then going to specify a sequence $(\ct_n)$ of 
finite subsets of the set $\ct$, on which we can base our
constructions.
The points in $\ct_n$ will consist of the centers of a minimal covering of
$\ct$  by closed 
balls of radius $\theta \rho^{n-1}$. 
It should be remarked that each set $\ct_n$ is not uniquely determined
and it is also possible that two 
points $s_m, t_n$ from the sets $\ct_m$ and $\ct_n$ for $m \neq n $
may be equal. This will not cause any trouble as one can see below.
 Later the numbers $\delta, \theta
\text{ and } \rho $  will determine which pairs of centers from the
union $\cup \ct_n$ we will
use for the construction of the two-point spectral triples, which will
be the summands in the spectral triple $ST(\delta)$.

\smallskip
\noindent
Given $\delta $, then the  smaller
 it is, the more accurate the metric $d_{\delta} $ describes the
 original metric. A better approximation to the given $d$ can only be
 achieved at a cost of having more points in the model.
 The exact content of this statement is coded in the integer
 $l(\delta)$ which we define below. We use the term {\em interaction
   length } for $l(\delta)$. It is not possible to see the meaning of
 this number right away, but it follows - hopefully - as the proof
 proceeds. The diameter of the space is given as diam$(\ct) \, = \,
 \max\{d(s,t)\,|\, s,t \in \ct \, \}$ and its logarithm 
 plays a role for the size
 of $l(\delta)$ so we define an integer $k_0$ - not necessarily non
 negative -  by
\begin{equation} \label{k0}
k_0 \, \in \, \bz,\quad \theta \rho^{k_0+1} \, < \, \text{diam}(\ct) \, \leq \theta
\rho^{k_0}.
\end{equation}
\noindent
Then we can define the interaction length $l\, := \, l(\delta)$ by the
conditions:

\begin{align} \label{l}
\text{If } \frac{4}{1-\rho}\, &<\, \delta \text{ then } l \, = \,
\max\{0, -k_0\} \\
\text{If } \frac{4}{1-\rho}\, &\geq \, \delta \text{ then } l
\, \in \bn \text{ is the least natural number s. t. }\\ \quad \, l &\geq -k_0
\text{ and } 
\frac{4\rho^l}{1-\rho}\, <\, \delta. 
\end{align}

\noindent
Remark that even though $k_0$ may be negative  the inequalities above
imply that 
\begin{equation} \label{kplusl}
k_0+l \geq 0 \text{ and } 
\frac{4\rho^l}{1-\rho}\, <\, \delta,  
\end{equation}
which will be useful later.
The two-point modules which will go into the
construction of the spectral triple $ST(\delta)$ can then be
determined.  As index set
$\cj$ we will
consider  all pairs of points $\{x,y\}$ from  $\ct$ such that     
\begin{equation} \label{indexpair} \end{equation} 
\begin{itemize}
\item[(i)] $\exists n \in \bn: x \in \ct_n \text{ and } y \in  \ct_n
  \cup \ct_{n+1}.$
\item[(ii)] $ x \neq y. $
\item[(iii)] If $ y \in \ct_n $ then $d(x,y) \leq
  \big(2+\rho^{-(l+1)}\big)\theta \rho^{n-1} $ 
\item[(iv)]  If $y \in \ct_{n+1}$ then $ d(x,y) \leq \big(1+\rho\big)\theta
  \rho^{n-1}.$
\end{itemize}
\noindent
Now our first task is to show that the direct sum of all these
two-point spectral triples will give a spectral triple. 
The first part of the proof of this 
 is done as in the proof of Theorem \ref{extripd}, but we remind the
reader that we are now using the standard two-point spectral triples
as defined in Definition \ref{twopttrip}.   
We will let $H$ denote the Hilbert space sum of all the $H_{x,y} $
for all the pairs $\{x,y\} \in \cj$.  This sum is countable  
and we may then define a
self-adjoint operator  $D$ on $H$ as the closure of the sum of all the
operators $D_{x,y} \text{ on } H_{x,y}$ over  $\{x,y\} \in \cj$. We will let $\pi$ denote the
representation of C$(\ct)$ on $H$ which 
 is equal to the sum of the representations $\pi_{x,y}$ on the spaces
$H_{x,y}, \text{ for } \{x,y\} \in \cj.$ 

\smallskip \noindent
As in the  proof of Theorem \ref{extripd} we  define the
set $\cn$, by 
\begin{displaymath}
\cn \,: =\, \{f \in
\text{C}(\ct) \, | \, \exists y \in \ct \, \forall x \in \ct\,:\,  f(x)=
d(x,y) \, \}.
\end{displaymath}
and you will find, as before,  that the algebra generated by  $\cn $
and the unit $I$ is uniformly dense in C$(\ct)$,  
and for any element, $f$ in this algebra, 
  the commutator  $[D,\pi(f)]$ is
bounded.
  The first condition in Definition \ref{spectrip}
 is then fulfilled. The second
condition, which asks for the compactness of the operator
$(I+D^2)^{-1}$ has to be considered in more details in order to be
verified. In order to prove that $(I+D^2)^{-1} $ is compact we just
have to show that for any positive number $r$ the self-adjoint 
operator $D$ has
only finitely many eigenvalues in the interval $[-r,r]$.
Let $\{x,y\}$ be a pair from $\cj$ then the eigenvalues, i. e. the
spectrum $\sigma(D_{x,y})$ of
$D_{x,y} $ is the set $\left \{ -d(x,y)^{-1}, d(x,y)^{-1} \right \}$, so we get that 
\begin{displaymath}
\text{If } x, y \in \ct_n,\,  \{x,y\} \in \cj  \text{ and }
\lambda \in \sigma(D_{x,y}) \text{ then }
 |\lambda| \geq  \frac{\rho^{1-n}}{(2+\rho^{-(l+1)})\theta}.
\end{displaymath}
\begin{displaymath}
\text{If } x  \in \ct_n,  y \in \ct_{n+1},\, \{x,y\} \in \cj  \text{ and }
\lambda \in \sigma(D_{x,y}) \text{ then }
 |\lambda| \geq \frac{\rho^{1-n}}{(1+\rho)\theta}.
\end{displaymath}
In both of the cases we see, that for any given $r> 0$ there will be
an $n_0$ such that only the finitely many elements $\{x, y \} \in \cj$
for which $x \in \ct_n$ for some $n \leq n_0$ can yield operators
$D_{x,y}$ with eigenvalues of absolute value less than or equal to
$r$. 

\medskip
\noindent
We will now turn to the properties of the metric 
induced by $ST(\delta)$. The first inequality claimed,
$d(x,y) \leq d_{\delta}(x,y) $, is quite general and it is proved exactly
as in  the proof of Theorem \ref{extripd}

\medskip
\noindent
We will then show the inequality $d_{\delta}(s,t) \leq (1+\delta)d(s,t)$.
Let $s,t$ be given from $\ct$, then we are first going to determine
the scale of the argument for this particular pair $(s,t)$,
  and this is done by finding the unique
integer $k$, it may be negative, such that 
\begin{equation} \label{k}
\theta \rho^k \, < \, d(s,t) \leq \theta \rho^{k-1}.
\end{equation}
\noindent 
By definition of $k_0$ we find that $k > k_0$ and by equation
(\ref{kplusl}) we then get 
\begin{equation} \label{k+l}
k + l \geq 1
\end{equation} 
\noindent
For any natural number $n$ we can by equation (\ref{k+l}) choose 
points $a_n, b_n \in  \ct_{k+l+n}$ such that
\begin{equation} 
d(s, a_n )\, \leq \, \theta \rho^{(k+l+n-1)} \text{ and } d(t, b_n)
\, \leq \, \theta \rho^{(k+l+n-1)} .
\end{equation}
Based on this we get via the triangle inequality 
\begin{align} \label{n,n+1}
\forall n \in \bn:\,d(a_{n+1},a_n)\,& \leq \, \theta
\rho^{(k+l+n)}+\theta \rho^{(k+l+n-1)} \text{ and } \\  d(b_{n+1},
b_n) \,& \leq \, \theta \rho^{(k+l+n )}+\theta \rho^{(k+l+n-1)}.
\end{align}
\begin{equation} \label{1,1}
d(a_1, b_1 )\, \leq \, d(s,t) + 2\theta \rho^{k+l}.
\end{equation}
from the last one and the definition of $k$ we get
\begin{equation} \label{x,y}
d(a_1, b_1 )\, \leq \, \theta\rho^{k-1} + 2\theta \rho^{k+l} =
\theta \rho^{k+l}\left ( 2 + \rho^{-(l+1)} \right ).
\end{equation}
Since  both $a_1$ and $b_1$ belongs to $\ct_{(k+l+1)} $ the last
inequality implies that the pair $\{a_1, b_1\} $ is in $\cj$ and 
the two-point spectral triple $ST_{a_1,b_1} $ is a summand in the spectral
triple we consider.

\smallskip
\noindent
Next we are going to show that the pairs $\{a_n,a_{n+1}\}, \{b_n,
b_{n+1}\} $ also belong to $\cj$. Let us look at $\{a_n, a_{n+1}\}$
then $a_n \in \ct_{k+l+n} \text{ and }  a_{n+1} \in \ct_{k+l+n+1}$, so 
we must check their distance. According to (4)  we get 
\begin{equation}
d(a_{n+1},a_n) \, \leq \, \theta
\rho^{(k+l+n)}+\theta \rho^{(k+l+n-1)} \, = \,\theta \rho^{(k+l+n-1)
}(1+\rho),
\end{equation} 
and it follows that these pairs also  go into the formation of
$ST(\delta)$, and analogously the  pairs $ \{b_n,
b_{n+1}\} $ belong to $\cj$ too.

\noindent
We can  then collect the estimates. The basic idea is that we  jump
from $s$ to an $a_n$, which is  nearby, and then continues to jump
from $a_i $ to 
$a_{i-1}$ until we reach $a_1$. From here we jump to $b_1 $ and
continues the jumping up to $b_n$ from where there is only a short
distance to $t$. All the jumps except the first and the last are
controlled by $ST(\delta)$. 
The end jumps can be made of arbitrary small
importance by continuity arguments. To be precise let us choose a
function $f \in \cam$, i. e. a continuous 
function such that $\|\, [D, \pi(f)]\,
\| \leq 1$ and try to get an upper estimate for $|f(s)-
f(t)|$. For a given $\eps > 0 $ there exists an $n \in \bn$ such that 
$|f(s)- f(a_n)| \leq \eps/2$ and $|f(t)- f(b_n)| \leq \eps/2$.
By construction the function $f \in \cam$ has the property that for
any pair $\{x,y\} \in \cj $ we have $|f(x) - f(y)| \leq d(x,y)$, so 
\begin{align*}
&|f(s) - f(t)| \leq \,  |f(a_n) - f(b_n)| + \eps \\
&\leq \,  |f(a_1) - f(b_1)| + \sum_{i=1}^{n-1}\big(|f(a_i) -f(a_{i+1})|
+ |f(b_i) -f(b_{i+1})|\big)+ \eps \\
& \leq \,  d(a_1, b_1) + \sum_{i=1}^{n-1}\big(d(a_i ,a_{i+1})
+ d(b_i,b_{i+1})\big)+ \eps \\
&\leq \,  d(s,t) + 4\theta\sum_{i=1}^{\infty}\rho^{(k+l-1+i)} + \eps\\
&= \, d(s,t)+ \frac{4\theta \rho^{k+l}}{1-\rho} + \eps \text{ and by
  } (\ref{k})\\
&\leq \, d(s,t) \left ( 1+ \frac{4 \rho^{l}}{1-\rho} \right ) + \eps \text{
  which by } (\ref{kplusl})\\
&\leq \, d(s,t)(1+ \delta) + \eps.
\end{align*}
Hence $d_{\delta}(x,y)\,  \leq \, (1+\delta)d(s,t)$.

\smallskip \noindent
We will then prove that $d_{\delta}$ generates the w*-topology on the
state space $\cs(\text{C}(\ct))$. This is also done as in the proof of
Theorem \ref{extripd} so we have to show that the  set $\cam_v$, which
is defined as in the proof of Theorem \ref{extripd} in the relations
(\ref{msepst}) and (\ref{camv}) is relatively compact and separates the states
of $\text{C}(\ct)$. The proof here is practically the same as the one
given at the end of the proof of Theorem \ref{extripd}. The only difference
being that we get the equicontinuity and the boundedness of $\cam_v$
from the following inequality.

\begin{displaymath}
\forall f \in \cam_v \,\forall s, t \in \ct:\quad
|f(s) - f(t) | \, \leq \, (1+\delta)d(s,t).
\end{displaymath}

\medskip
\noindent
We will now turn to the summability questions. We think that the
results obtainable will turn out to be much more precise in the future,
 but we have gotten into a lot of combinatorial problems when we tried to
sharpen our  results on the connection between finite  summability of
$ST(\delta)$ and finiteness of the upper Minkowski dimension of $\ct$. 

\noindent
Let us suppose that the upper Minkowski dimension,
$\overline{\text{dim}_{\text{M}}}(\ct)$  is
finite, and recall that for $n \in \bn$ the set $\ct_n$ consists of
the centers of a minimal covering of $\ct$ by closed balls of radius
$\theta\rho^{n-1}$. We will let $|\ct_n|$ denote the number of points
in $\ct_n$ then we get by Definition \ref{upmindim} 
\begin{displaymath}
\overline{\text{dim}_{\text{M}}}(\ct)\, := \, \underset{n\to
  \infty}{\limsup}\frac{\log(|\ct_n|)}{-\log(\theta\rho^{n-1})} < \infty.
\end{displaymath} 
Let $\mu $ be a real number such that $\mu >
\overline{\text{dim}_{\text{M}}}(\ct)$ then there must exist a natural
number $n_0$ such that

\begin{displaymath}
\forall n > n_0: \quad
\frac{\log(|\ct_n|)}{-\log(\theta\rho^{n-1})}\, < \,\mu, 
\end{displaymath}
and then there must be a natural number $n_1 \geq n_0$ such that
\begin{displaymath}
\forall n > n_1: \quad  \frac{\log(|\ct_n|)}{-\log(\rho^n)}\, < \, \mu,
\end{displaymath}
which implies that
\begin{displaymath}
\forall n > n_1: \quad  |\ct_n| \, < \,\rho^{-\mu n}.
\end{displaymath}
We can now make estimates of the value of tr$\big((I+D^2)^{-s/2}\big)$ and we
will show that it is finite for any real $s >
2\overline{\text{dim}_{\text{M}}}(\ct)$. Suppose such an $s$ is given
and $\mu $ is chosen such that $2\overline{\text{dim}_{\text{M}}}(\ct)
< 2\mu < s$. First we remark that for $\lambda \neq 0$  we have the inequality 
\begin{displaymath} 
\frac{1}{\sqrt{1+|\lambda|^2}}\, \leq \, |\lambda|^{-1}.
\end{displaymath} 
When we have to estimate the trace  $\text{tr}\big((I+D^2)^{-s/2}\big)$ we must go back
to the definition of $D$, and we find that the index set $\cj,$ over
which $D$ is formed, is grouped into disjoint subsets $\cj_n, \, n\in
\bn$ such that 
\begin{displaymath}
\cj_n \, := \, \{\, \{x, y \} \in \cj \, | \, x \in \ct_n \text{ and }
y \in \ct_n \cup \ct_{n+1} \, \}.
\end{displaymath}
Having this and the inequality just above we get
\begin{displaymath} 
\text{tr}\big((I+D^2)^{-s/2}\big) \, \leq \, \sum_{n=1}^{\infty}\left ( \,
\sum_{\{x,y\} \in \cj_n} 2d(x,y)^s\,\right ).
\end{displaymath}
In this sum we recall that there are upper limits for $d(x,y) $ which
were used in the definition of the spectral triple. By inspection
of the definition it turns out that we have the following inequalities 
\begin{displaymath} 
 \exists c> 0\forall n \in \bn \forall \{x,y\} \in \cj_n\,: \quad d(x,y) \leq c\rho^n.
\end{displaymath} 
The number, say $|\cj_n|$, of pairs in $\cj_n$
 is not known, because  it depends of the
local structure  of the metric $d$. Without further knowledge we can
only get some very rough estimates of the size of $|\cj_n|$, but we
can remark that for a nice subset of the space $\br^d$ such as the
unit cube, one can do much better. The general estimate we can get is simply
the largest possible numbers of pairs, so we get
\begin{displaymath}
\forall n \in \bn: \quad |\cj_n| \leq \frac{|\ct_n|^2}{2} +
|\ct_n||\ct_{n+1}|.
\end{displaymath} 
This can then be combined with the estimates $|\ct_n| < \rho^{-\mu n}$,
for $n > n_1$, so we obtain
\begin{align*}
&\text{tr}\big((I+D^2)^{-s/2}\big)\\ & \leq  \sum_{n=1}^{n_1} \left ( \,
\sum_{\{x,y\} \in \cj_n} 2d(x,y)^s\,\right ) + \sum_{n>n_1} \left ( \,
\left ( c\rho^n \right ) ^s \left ( \frac{\rho^{-2\mu n}}{2} + \rho^{-\mu(2n+1)} \right ) \right ).
\end{align*}
This means that there exist positive reals $A$  and $B$ such that 
\begin{displaymath}
\text{tr}\big((I+D^2)^{-s/2}\big)  \leq  A + B\sum_{n>n_1}\rho^{(s-2\mu)n} <
\infty,
\end{displaymath}
and the module is finitely summable for any $s >
2\overline{\text{dim}_{\text{M}}}(\ct)$.  

\medskip
\noindent
Let us now suppose that the module is finitely summable for some $s >
0$  and the space $\ct$ is connected. By the construction of the
spectral triple (\ref{indexpair}) a pair of centers $\{x,y\}$ in
$\ct_n$  belongs to $ \cj_n$ if  
$d(x,y) \leq  \left ( 2+\rho^{-(l+1)} \right ) \theta\rho^{n-1}. $ This means, among
other things,  that there exists a natural number $n_2$ such that for
any natural number $n > n_2$ and any pair $\{x,y\} \in \ct_n$ such
that $x,y \in \cj_n$ then we have

\begin{displaymath} 
d(x,y) \, \leq \, 1, \quad \text{and} \quad 2\theta\rho^{n-1} \, \leq \,
\text{diam}(\ct)/3. 
\end{displaymath}
 The inequality $d(x,y) \leq 1 $ implies $\left ( 1+d(x,y)^{-2} \right ) ^{-s/2}  \, \geq
 2^{-s/2}d(x,y)^s, $ so the assumption of finite summability of the
 module and these considerations imply

\begin{align}
\infty \,>&\, \text{tr}\big((I+D^2)^{-s/2}\big) \\ =&\, \sum_{n=1}^{\infty}\left ( \,
\sum_{\{x,y\} \in \cj_n} 2 \left ( 1+d(x,y)^{-2} \right ) ^{-s/2}\,\right )\\
\geq & \, \sum_{n>n_2}\left ( \,
\sum_{\{x,y\} \in \cj_n, x,y \in \ct_n} 2^{(1-s/2)}d(x,y)^s\,\right ).  \label{lboundtr} 
\end{align}
We are now going to show that for each $n > n_2$ there are at least
$|\ct_n|/2$ number of pairs $\{x,y\}$ in $\cj_n$ such that $x$ and $y$
both are in $\ct_n$ and 
$d(x,y)
\geq \theta\rho^{n-1}$. We do this by showing that for  any $x \in
\ct_n$ there exists at least one $y \in \ct_n$ such that $\{x,y\} \in
\cj_n$ and $d(x,y) \geq \theta\rho^{n-1}$. Let  $x \in \ct_n$ then
there must exist a $t \in \ct $ such that $d(x,t) \geq
\frac{\text{diam}(\ct)}{3}$. By assumptions $\ct$ is connected and
for $n>n_2,\, 2\theta\rho^{n-1}\, \leq \,
\frac{\text{diam}(\ct)}{3}$, hence there must be a $u \in \ct$
such that $d(x,u) = 2\theta\rho^{n-1} $. Finally we choose $y \in \ct_n
$ such that $d(y,u) \leq \theta\rho^{n-1}$ and we get 
\begin{displaymath}
 d(x,y)\, \leq  \, d(x,u) + d(u,y) \, \leq
3\theta\rho^{n-1} <  \left ( 2 + \rho^{-(l+1)} \right )\theta\rho^{n-1},
\end{displaymath}
\noindent
so by (\ref{indexpair}) the pair $\{x,y\}$ belongs to the index set
$\cj_n$, if $x \neq y$.
 On the other hand, if we  use the triangle
inequality the other way around, we get 
\begin{displaymath}
d(x,y) \, \geq \, d(x,u) - d(u,y) \, \geq \, \theta\rho^{n-1} > 0,
\text{ so } x \neq y.
\end{displaymath} 
\noindent
We have then seen that any point $x$ in $\ct_n$ belongs to a pair
$\{x,y\}$ in $\cj_n$ with $y$ in $\ct_n$ and  $d(x,y) \geq
\theta\rho^{n-1}$, so  we may continue our
lower bound estimates (\ref{lboundtr})  on the 
trace  of $(I+D^2)^{-s/2}$ in order to get the following result
\begin{displaymath}
 \sum_{n > n_2} |\ct_n|(\rho^{s})^n \, < \, \infty.
\end{displaymath}
In particular this means that there exists a natural number  $n_3 > n_2$
such that for any $n > n_3$ 
\begin{align*}
 |\ct_n|(\rho^{s})^n \, <& \, 1, \text{ and } \\
\log |\ct_n| \, < & \, -s \log \rho^n \\
\frac{\log |\ct_n|}{- \log \rho^n} \, < & \, s, \text{ and for n
  suff. large} \\
\frac{\log |\ct_n|}{- \log(\theta \rho^{n-1})} \, < & \, s.
\end{align*}
\noindent
Since the set $\ct_n$ is a set of centers for a minimal covering of
$\ct$ by balls of radius $\theta\rho^{n-1}$, the last inequality shows
that the upper Minkowski dimension of $\ct$ is at most s, and the
theorem follows.
\end{proof}

\noindent
\begin{rem}
The result raises a number of questions. First of all, is it possible
to construct {\em discrete } spectral triples which are more accurate
with respect to the calculation of the Minkowski dimension ? 
 We think that the answer
definitely is yes, and for all the known examples we have thought of,
there are ad hoc constructions which are very precise.
The works  by Guido and Isola \cite{GI1, GI2, GI3}  offer spectral
triples which works for a large collection of fractals in
$\br^n$. Some of the problems we have faced when we have tried to get
more precise results is - in the language of the proof just above - 
to determine the number of pairs $\{x, y\} \in \cj_n$ which contains a
given $x$. In known examples there exists a constant $c $ such that
this number is at most $c$ for any value of $n$. In our proof of
Theorem \ref{apptripd} the
estimate is 3/2$|\ct_n|^2$ which is far away from being universally
bounded. To illustrate the exact content of these cryptical remarks we
will compute an example based on the unit interval. A similar result
can be obtained for any unit cube in $\br^d$, but the counting of the 
multiplicities of the eigenvalues is much more complicated in higher
dimensions, so in order to avoid huge 
combinatorial arguments we will just compute the Dirac operator $D$  for
the unit interval based on the construction given in the proof of 
Theorem \ref{apptripd}.  We then show that this
spectral triple behaves very nicely.  The only un-excepted result is
that for this module the so called Weyl's asymptotic formula for the
growth of the eigenvalues does not hold completely.
 The reason for this is that this Dirac operator has gaps between the
 eigenvalues which grow exponentially and also the multiplicities of
 the eigenvalues grow exponentially. You can see the details just
 below. 
\end{rem}

\begin{exmp}
We consider the unit interval $[0,1]$ with the usual metric. In the
construction from the proof of Theorem \ref{apptripd} we let $\delta = 9,
\, \theta = 1  $ 
and $\rho = 1/2. $ Then for the spectral triple $ST(9) = \big(
{\rm C}([0,1]), H, D \big)$ constructed by
the algorithm of Theorem \ref{apptripd} we get:
\begin{itemize}
\item[(a)] The metric induced by the spectral triple is the usual
  one. 
\item[(b)] The module is finitely summable for $s>1$, but not for $s=1$.
\item[(c)]  Let for $\Lambda > 0, N(\Lambda)$ denote the number of
  eigenvalues of $D$ of absolute value at most $\Lambda$ then
$$  10 \, < \, \underset{\Lambda \to \infty } {\liminf}
 \frac{N(\Lambda)}{\Lambda} \, \leq \, 13 \, < \, 17 \,
 \leq \,  \underset{\Lambda \to
   \infty } {\limsup } \frac{N(\Lambda)}{\Lambda}\, < \, 20. $$

\item[(d)] For $\Lambda > 0$ let $P_{\Lambda}$ denote the orthogonal
  projection onto the subspace of $H$ spanned by the eigenvectors of
  $D$ with non zero eigenvalues of numerical value at most $\Lambda.$
  For any   any continuous function $f$ on $[0,1]$ we get  
$$\underset{\Lambda \to \infty }{\lim} \frac{1}{\log(\Lambda)}{\rm
  tr}\left (  |D|^{-1}P_{\Lambda}\pi(f)\right )  =  \frac{10}{\log 2}\int_0^1f(x){\rm d}x.$$
\item[(e)] For any ultrafilter $\omega $ on $\bn$ and any continuous
  function $f$ on $[0,1]$ the Dixmier trace for this spectral
  triple is given by $${\rm tr}_{\omega}\left ( \pi(f)|D|^{-1} \right )= \frac{10}{\log 2}\int_0^1f(x){\rm d}x.$$  
\end{itemize}
A different choice of $\delta,
\theta $ and $\rho$ would give a similar result, except for the
constants $\log2, 10, 13, 17, 20.$  

\end{exmp}
\noindent
{\bf Details of the example}
\noindent
The hole construction is based on the sets $\ct_n$, which for a given
$n$ is a set of centers for a minimal covering of the interval $[0,1]
$ by balls of radius $\theta \rho^{n-1} = 2^{1-n}$. So $\ct_1 = \ct_2 =
\{2^{-1}\}, $ and for 
$$ n\,  > \, 2 \quad \ct_n \, = \, \left \{ \, (2j + 1) 2^{1-n}\, | \, 0 \leq
j \leq 2^{n-2} -1 \,\right \}.$$
In particular the number of points $|\ct_n| = 2^{n-2} $ for $n \geq
2.$
\smallskip \noindent
We can then prove the statement (a) without too much trouble. Let $u <
v$ be points in $[0,1]$ and let $\eps > 0$ be given. Let 
 $d_{\d}$ denote the metric induced by the spectral triple. We
then know from the theorem, 
that for any pair of points $u, v$ from the unit interval we
have $d_{\d}(u,v) \leq |u-v|.$ Let then $f$ be a continuous function
such that  $\|\, [D,\pi(f)]\,\| \leq 1,  $ and let $\eps > 0$ be
given. Then by the continuity of $f$ we can find 
a natural number $n$ and points $x_j <
x_k$ from $\ct_n$ such that $|f(u)-f(x_j)| < \eps/2$ and $|f(v) - f(x_k)|
< \eps/2$. By the choice of $f$ we have $$|f(x_k) - f(x_j)| \leq
\sum_{m=j}^{k-1} |f(x_{m+1}) - f(x_m| \leq  \sum_{m=j}^{k-1} |x_{m+1}
- x_m| = x_k - x_j.$$ \noindent Hence $|f(v) - f(u)| \leq |u-v| +
\eps$ and claim (a) follows.

\smallskip
\noindent
The next task is then for a given natural number $n$ and each point $x
$ in the   set $ \ct_n$ to
study  the one dimensional representation, say
$\pi_x$, given as  
C$([0,1]) \ni f \to f(x) $ as a summand in  representation $\pi$ of the
spectral triple. Suppose that a pair $\{x,y\} $ belongs to the index
set $\cj$ then both $\pi_x $ and $\pi_y $ will be counted once {\em
  more } in $\pi$ so the multiplicity of $\pi_x$ in $\pi $ is exactly
the number of pairs $\{u,v \} \in \cj$ which contains $x$. To find
that number we will first determine the numbers $k_0$ and $l$ used in
the proof of Theorem \ref{apptripd}. According to the relation
(\ref{k0}) we get $k_0= 0$ and by (\ref{l}) we get $l=0$ too.
We then have to  look at the rules (\ref{indexpair}) which defines the
pairs in $\cj$, and one gets 
for an $ n \geq 5 $  and an  $x_j = (2j+1)2^{1-n}$ in $\ct_n$ that we
have to look for points $v_i = (2i+1)2^{2-n} \in \ct_{n-1}, y_k =
(2k+1)2^{1-n} \in \ct_n $ and $z_l = (2l+1)2^{-n} \in \ct_{n+1}$ such
that the pairs $\{v_i,x_j\}, \, \{x_j,y_k\},\, \{x_j, z_l\}$ 
in which $x_j$ appears, are   determined  by the following relations 
\begin{align*}
|v_i- x_j| \quad \leq &\quad \frac{3}{2}\cdot 2^{2-n} \quad = \quad  6 \cdot 2^{-n}\\
|x_j - y_k| \quad \leq &\quad 4\cdot 2^{1-n}         \quad  =  \quad 8\cdot 2^{-n}\\
|x_j - z_l| \quad \leq & \quad \frac{3}{2}\cdot 2^{1-n}   \quad  =    \quad 3\cdot 2^{-n}.
\end{align*}
\noindent
There is a pattern for the computation of this, but it only holds for
an $x_j$ which is not very near $0$ or $1$. A detailed investigation
shows that the regular pattern is valid for all $x_j$ for which $j \in
\{2,  \dots, 2^{n-2}- 3\}$, so it is only 4 points from each of the
sets $\ct_n$ which we will forget in the computations to come. Further
it is not difficult to realize 
that the properties of the module is not dependent on
any finite number of summands. Any finite number of two-point modules
may be omitted without disturbing
the properties of the module with respect to the metric or any of the
asymptotic properties which we want to prove now. The reason for this
is that any sort of property is always better described for a higher
$n.$
We may, and will,
then assume that the module only starts by the summation of all the
modules which involves points from $\ct_n$ for $n \geq 5.$ This means
that some points from $\ct_4$ will be introduced, as they are part of
some two-point modules with points from $\ct_5.$ Let now $x_j =
(2j+1)2^{1-n} $ be a point from $\ct_n$ for an $n \geq 5$ and a $j$
such that $2 \leq j \leq 2^{n-2} - 3.$ We then get the following table
for the number of points.
\begin{align*}
& | \{ v_i \in \ct_{n-1}\, |  \, |v_i- x_j| \leq \frac{3}{2}\cdot 2^{2-n}\,  \} | 
=2 \\ & \text{ and the distances between } x_j \text{ and } v_i 
\\& = \left ( 2\cdot 2^{-n}, 
6\cdot 2^{-n} \right ).\\ 
& | \{ y_k \in \ct_{n}\, | \, |y_k- x_j| \leq 4\cdot 2^{1-n}\, \} |  = 4 \\
& \text{ and the distances  between } x_j \text{ and } y_k
\text{ and repeated if necessary }\\  & =
\left ( 4\cdot 2^{-n},4\cdot 2^{-n},8\cdot 2^{-n},8\cdot 2^{-n}\right
).\end{align*} \begin{align*} 
& | \{ z_l \in \ct_{n+1}\, | \, |z_l- x_j| \leq
    \frac{3}{2}\cdot 2^{1-n}\, \} | = 4 \\ & \text{ and the distances  between } x_j \text{ and } z_l
\text{ and repeated if necessary }\\ & = \left ( 2^{-n},
2^{-n}, 3\cdot 2^{-n}, 3\cdot 2^{-n} \right ).
\end{align*} 
\noindent
With respect to a  point $x_j$ in $\ct_n$ we then get that the multiplicity
of $\pi_x$ is at most $10$, and it is $10$ for all the points in
$\ct_n$ except for the 4 points nearest the boundary $\{0, 1\}$.
 
\smallskip
\noindent
We then turn to the absolute value $|D|$ of the Dirac operator. Let us
then look at a pair $\{x,y\} \in \cj$, then the Definition \ref{twopttrip}
shows that $|D_{x,y}| = |x-y|^{-1}I_{H_{x,y}}$ so a  unit vector in
the Hilbert space $H_{x,y} $
 corresponding to the representation $\pi_{x}$ is
an eigenvector for $|D| $  corresponding to the eigenvalue
$|x-y|^{-1}$.  This means that for an $x_j = (2j+1)2^{1-n},$ with $2
\leq j \leq 2^{n-2} -3 $ the $10$ eigenvalues  for the operator  $|D|$
corresponding to the summands of $\pi_{x_j}$ in $\pi$ will be; when
listed by multiplicity and increasingly,

  $$ \left \{\, 2^{n-3},\,\, 2^{n-3},\,\,  \frac{ 2^{n-1}}{3}, \,\, 2^{n-2}, \,\, 2^{n-2}, \,\,  \frac{2^{n}}{3}, \,\,
\frac{2^{n}}{3}, \,\,   2^{n-1}, \,\,
2^{n}, \,\,   2^{n} \, \right \}.$$
  
\medskip
\noindent
We can  now begin to compute the multiplicities of the eigenvalues
for  $|D|.$ We know that all  of the
 $2^{n-2} $ points $x_j$ in the set $\ct_n$,  except $4$,  behave
 alike. For these 4
 points the  set of eigenvalues for $|D|$ are the same, but with lower
 multiplicities. In the computations to come it makes no difference, when
 estimating limits, to neglect this irregularity, so we will just
 assume that the 4 outmost points of each set $\ct_n$ 
have the property that the absolute value of the Dirac
 operator follow the same eigenvalue pattern here, as for the regular points.
From the list above we find that
the only possible eigenvalues for $|D|$ on a vector corresponding to
$\pi_{x_j}$ are of the forms $  \frac{ 2^k}{3}$ and $
2^l$. Further we see from the list, that the eigenvalue $2^n$ can only
be associated with points from $\ct_n, \ct_{n+1}, \ct_{n+2}$ and
$\ct_{n+3}$ and we get, since the number of points in a set $|\ct_k|$ is
$2^{k-2}$, that the multiplicity of  the eigenvalue $2^n$  will be
nearly 
$$2\cdot 2^{n-2}\quad + \quad 1\cdot 2^{n-1}\quad + \quad 2\cdot 2^{n}\quad +
\quad 2\cdot 2^{n+1}\quad = \quad 7\cdot 2^n$$ 
  
\noindent 
Similarly the eigenvalue $  \frac{ 2^n}{3}$, will be generated by
points in $\ct_n$ and  $\ct_{n+1}$ and the multiplicity of $  \frac{
  2^n}{3}$ for $|D|$ is approximately

$$2\cdot 2^{n-2}\quad + \quad 1\cdot 2^{n-1}\quad = \quad 2^n.$$ 
\noindent
We then get 
\begin{align*}
\text{tr}(|D|^{-s} ) \, =& \, \sum_{n=5}^{\infty}\left ( 7\cdot 
2^{(n-ns)} + 3^s\cdot 2^{(n-ns)} + E(n,s)2^{-ns}\right ) \text{ such that }\\
|E(n,s)| \,  \leq& \, 40\cdot 3^s.
\end{align*}

\noindent It then follows that the module is summable for $s > 1$ and
not summable for $s = 1.$ The item (b) of the example is then proved,
and we will begin to estimate $N(\Lambda)$, the number of eigenvalues
for $|D|$, dominated by $\Lambda$. We can base the computation on the
counting of the multiplicities which we just performed above. 
Let $\lfloor \frac{\log \Lambda}{\log 2} \rfloor$ denote the largest
integer smaller than or equal to   $\frac{\log \Lambda}{\log 2},$ then
the number of eigenvalues $N(\Lambda)$ for $|D|$, which are dominated
by $\Lambda $ is nearly
\begin{equation} \label{N}
\sum_{n=5}^{\lfloor \frac{\log \Lambda}{\log 2} \rfloor}7\cdot 2^n \quad
+ \quad \sum_{n=5}^{\lfloor \frac{\log 3\Lambda}{\log 2} \rfloor}2^n
\end{equation}
\noindent
and we get   
\begin{align*}
 N(\Lambda)\quad =& \quad 7 \cdot 2^{\lfloor \frac{\log \Lambda}{\log 2} \rfloor
   +1}\quad  + \quad 2^{\lfloor \frac{\log 3\Lambda}{\log 2} \rfloor
   +1} \quad + \quad R(\Lambda) \text{ such that } \\ 
|R(\Lambda)| \quad <& \quad  40\frac{\log \Lambda +3}{\log2} + 100
\quad <  \quad 80
\log(\Lambda) + 500.
\end{align*}
\noindent
We then get 
$$ N(\Lambda)\,  = \, \Lambda \left ( 7\cdot 2^{\left ( 1 - \left ( \frac{\log
    \Lambda}{\log 2}  
 -  \lfloor \frac{\log \Lambda}{\log 2}\rfloor \right ) \right )}\, + \, 3\cdot 2^{\left ( 1 -
 \left ( \frac{\log 3\Lambda} {\log 2} 
 -  \lfloor \frac{\log 3\Lambda}{\log 2}\rfloor \right ) \right )}\right ) \, + \,
R(\Lambda).$$
\noindent
The exponents for $2$ in the expression above are always between $0$
and $1$ and they will vary discontinuously between these values so
$\frac{N(\Lambda)}{\Lambda} $ has no limit for $\Lambda \to \infty$
and we get 
$$
 10 \, \leq \, \underset{\Lambda \to \infty } {\liminf}
 \frac{N(\Lambda)}{\Lambda} \, \leq \, 13 \, < \, 17 \,
 \leq \,  \underset{\Lambda \to
   \infty } {\limsup } \frac{N(\Lambda)}{\Lambda} \, \leq \, 20 $$   
\noindent   
This establishes the item (c) of the example. 
It seems likely that one can do better, but $\frac{\log 3}{\log 2} $
is not an integer so it is not possible to get down to $10 $ or up to
$20$.

\smallskip
\noindent
Let us then turn to item (d). Let $\eps > 0  \text{ and } f$  a
continuous complex
function on the interval $[0,1]$ be given.  Since $f$ is continuous  there
exists a natural number  $N_1$ such that
\begin{align} \label{int} 
\forall n \in \bn, \, n \geq N_1: \,& \left \vert \,\int_0^1 f(x)\text{d}x -
\sum_{j=0}^{2^{n-2}-1}f \left ( (2j+1)2^{1-n}\right ) 2^{2-n}\,\right \vert  \leq \frac{\eps}{120}
\\ \text{ and } & 40\|f\|2^{2-n} \leq  \frac{\eps}{12} \label{40f}
\end{align}
\noindent
Then we define for any natural  number $n \geq 5$ and any $j \in \{0,
\dots , 2^{n-2}-1\}$,   the point $x_j \in \ct_n$ as before by $x_j =
(2j+1)2^{1-n}$ and the inequality (\ref{int}) gives  
\begin{equation} \label{int-sum}
  \forall n \in \bn, \, n \geq N_1: \quad \left \vert \,\int_0^1 f(x)\text{d}x -
\sum_{j=0}^{2^{n-2}-1}f(x_j)2^{2-n}\,\right \vert  \leq \frac{\eps}{120}.
\end{equation}
Let now  $Q_{n,j}$ be the orthogonal  projection onto the span of
the vectors  in 
$H$ which have the property that each of them  
generates the one dimensional representation $\pi_{x_j}$
of C$([0,1])$. Let us then
 go back to the computations of the pairs which involved an $x_j, $
 and we find - again except for the 4 points which are closest to
 either $0$ or $1$ -  that for an $x_j$ from $\ct_n$, the multiplicity
 of $\pi_{x_j}$ in $\pi$ is $10$ and the corresponding $10 $ eigen
 values for $|D|^{-1} $ are  
$$ \left \{ 2^{-n}, \,\, 2^{-n}, \,\, 2^{1-n},  \,\, 3\cdot  2^{-n}, \,\,
3\cdot 2^{-n}, \,\, 2^{2-n}, \,\, 2^{2-n}, \,\, 3\cdot 2^{1-n}, \,\, 
2^{3-n}, \,\,  2^{3-n}\, \right \}$$ 
\noindent We can then for $2 \leq j \leq 2^{2n-1}-3$ compute 
$$
\text{tr}(Q_{n,j}|D|^{-1})  =
2^{-n}(1+ 1 +2+ 3 +3 +4 +4 + 6+8+8) = 10\cdot 2^{2-n},
$$
and for the 4 outmost points we get 
$$
2^{2-n} \, < \, \text{tr}(Q_{n,j}|D|^{-1})  \, < \,   10\cdot 2^{2-n},
$$
It then follows that for the projection $Q_n\, := \,
\underset{j=0}{\overset{2^{n-2}-1}{\sum}}Q_{n,j}$ and any continuous complex function
  $g$ on $[0,1]$ we have 
\begin{align} \label{Qnest}
\text{tr}(Q_n\pi(g)|D|^{-1})\, =& \,
10\sum_{j=0}^{2^{n-2}-1}g(x_j)2^{2-n} + I(n) \text{ such that } \\
|I(n)| \, \leq & \, 40 \|g\|2^{2-n} \text{ and we also get} \label{40g} \\
\left \vert \text{tr}\left ( Q_n\pi(g)|D|^{-1} \right ) \right \vert  \, \leq & \, 10\|g\|. \label{10g}
\end{align} 
\noindent
Hence for the given $f$ and a natural number $n \geq N_1 $ we may use 
(\ref{40f}, \ref{int-sum}, \ref{Qnest}, \ref{40g}) to obtain 
\begin{equation} \label{tr-int}
\left \vert \text{tr}\big(Q_n\pi(f )|D|^{-1}\big) \, - \,10 \int_0^1 f(x)\text{d}x \right \vert \,  \leq
\, \frac{\eps}{6}.
\end{equation}
\noindent
We can now introduce the positive real  $\Lambda$ and we will choose
it such that 
\begin{align} 
\Lambda \geq & \exp\left ( \frac{90(N_1 + 2)(\|f\|+1)}{\eps}\right ) \text{ and }
\Lambda \geq 2^5 = 32 \text{ so} \\
\frac{\eps}{90} \geq & \frac{(N_1 + 2)(\|f\|+1)}{\log \Lambda}. \label{eps}
\end{align}
 \noindent
Then we get for the given continuous function $f$ 
and for $M$ a natural number defined by $M:=\lfloor \frac{\log
  \Lambda}{\log 2}\rfloor$ that  
\begin{align}
\frac{1}{\log \Lambda}\text{tr}\big(P_{\Lambda}\pi(f)|D|^{-1}\big) \,
=& \,  \frac{1}{\log \Lambda}
\sum_{n=5}^{M}\text{tr}\big(Q_n\pi(f)|D|^{-1}\big) + J(\Lambda) \\
\label{31} \text{ such that }   
|J(\Lambda)| \, = & \,  \left \vert \frac{1}{\log \Lambda}\text{tr}\left (
 \left ( P_{\Lambda}- \sum_{n=5}^{M}Q_n \right )\pi(f)|D|^{-1}\right ) \right \vert \\
 \leq & \,   \left \vert \frac{\|f\|}{\log \Lambda}\text{tr}\left (
 \left ( P_{\Lambda}- \sum_{n=5}^{M}Q_n \right ) |D|^{-1}\right ) \right \vert , \text{  by
  (\ref{N})} \\
 \leq & \, \left \vert \frac{\|f\|}{\log \Lambda}\text{tr}\left (
 (Q_4 + Q_{M+1}+ Q_{M+2})|D|^{-1}\right ) \right \vert  \\
 \leq & \,  \frac{30\|f\|}{\log \Lambda} \\
 < &\frac{\eps}{90}. \label{J}  
\end{align}
We can then continue as 
\begin{align} \label{1/L}
& \frac{1}{\log \Lambda}
\sum_{n=5}^{M}\text{tr}\big(Q_n\pi(f)|D|^{-1}\big) \\ &= \,
\frac{1}{\log \Lambda}
\sum_{n=5}^{N_1}\text{tr}\big(Q_n\pi(f)|D|^{-1}\big) + \frac{1}{\log \Lambda}
\sum_{n=N_1+1}^{M}\text{tr}\big(Q_n\pi(f)|D|^{-1}\big) \\ 
& = \frac{M - N_1 -1}{\log \Lambda}10\int_0^1 f(x)\text{d}x + K(\Lambda)
 \text{ such that by  (\ref{10g}, \ref{tr-int}) }, \\ & |K(\Lambda)| \,
 \leq  \, \frac{10\|f\|N_1}{\log
  \Lambda} + \frac{M - N_1 -1}{\log \Lambda} \frac{\eps}{6} \leq
 \frac{\eps}{9}
+\frac{\eps}{6} < \frac{\eps}{3}. \label{39}
\end{align}
The final estimation then becomes 
\begin{align*}
& \left \vert \,\frac{1}{\log \Lambda}\text{tr}\left ( P_{\Lambda}\pi(f)|D|^{-1}\right ) -
\frac{10}{\log 2}\int_0^1 f(x)\text{d}x\,\right \vert  \, \, \text{ by
} (\ref{31} - \ref{J}) \\  
\leq & \left \vert \frac{1}{\log \Lambda}
\sum_{n=5}^{M}\text{tr}\left ( Q_n\pi(f)|D|^{-1}\right ) - \frac{10}{\log
  2}\int_0^1 f(x)\text{d}x\,\right \vert \, +\, \frac{\eps}{90} \,\,
\text{ by } (\ref{1/L} - \ref{39}) \\
\leq & \left \vert \,\frac{10(M - N_1 -1)}{\log \Lambda} -  \frac{10}{\log
  2}\,\right \vert  \, \left \vert \,\int_0^1 f(x)\text{d}x\, \right \vert  \,  +\, \frac{\eps}{3}\, +\,
\frac{\eps}{90}  \\
< & 10 \|f\| \frac{\log \Lambda - \log 2\left ( \lfloor\frac{\log
    \Lambda}{\log 2}\rfloor - N_1 - 1\right )}{ \log2 \log \Lambda} \, + \, \frac{\eps}{2}\\
\leq & 10 \|f\| \frac{\log \Lambda - \log 2\left ( \frac{\log
    \Lambda}{\log 2}-1 - N_1 - 1 \right )}{ \log2 \log \Lambda} + \frac{\eps}{2}  \\
= & 10 \|f\| \frac{N_1 + 2}{ \log \Lambda} \, + \,
\frac{\eps}{2}\text{ which by }  (\ref{eps}) \\
\leq & \frac{\eps}{9} \, +\, \frac{\eps}{2 } < \eps. 
\end{align*}
This establishes item (d) of the example and we will turn to look at
item (e), but when we have a convergence as described in item (d),
then it follows from the construction of the Dixmier trace \cite{Co2}
Chapter IV 2.$\beta$, pp 303 -- 308, that the limit obtained in (d) is
the Dixmier trace of $\pi(f)|D|^{-1}.$

\end{document}